
%
%
%
%
%

\input amstex
\loadbold
\documentstyle{amsppt}
\NoBlackBoxes

\pagewidth{32pc}
\pageheight{44pc}
\magnification=\magstep1

\def\QQ{\Bbb{Q}}

\def\PP{\Bbb{P}}
\def\OO{\Cal{O}}
\def\CC{\Bbb{C}}
\def\Spec{\operatorname{Spec}}
\def\Supp{\operatorname{Supp}}

\def\image{\operatorname{Im}}
\def\rank{\operatorname{rk}}
\def\codim{\operatorname{codim}}

\def\Mov{M}
\def\Height{\operatorname{Height}}
\def\achern#1#2{\widehat{c}_{#1}(#2)}
\def\rest#1#2{\left.{#1}\right\vert_{{#2}}}

\topmatter
\title
Faltings modular height and self-intersection of dualizing sheaf
\endtitle
\rightheadtext{}
\author Atsushi Moriwaki \endauthor
\leftheadtext{}
\address
Department of Mathematics, Faculty of Science,
Kyoto University, Kyoto, 606-01, Japan
\endaddress
\curraddr
Department of Mathematics, University of California,
Los Angeles, 405 Hilgard Avenue, Los Angeles, California 90024, USA
\endcurraddr
\email moriwaki\@math.ucla.edu \endemail
\date January, 1994 \enddate
\thanks
THIS IS A TENTATIVE VERSION.
\endthanks
\endtopmatter

\document

\head 0. Introduction \endhead

Let $K$ be a number field, $O_K$ the ring of integers of $K$ and
$X$ a stable curve over $O_K$ of genus $g \geq 2$.
In this note, we will prove
a strict inequality
$$
\frac{(\achern{1}{\omega_{X/S}, \Phi_{can}}^2)}{[K : \QQ]} >
\frac{4(g-1)}{g} \Height_{Fal}(J(X_K)),
$$
where $\omega_{X/S}$ is the dualizing sheaf of $X$ over $S = \Spec(O_K)$,
$\Phi_{can}$ is the canonical Hermitian metric of $\omega_{X/S}$
and $\Height_{Fal}(J(X_K))$ is the Faltings modular height of
the Jacobian of $X_K$ (cf. Corollary~2.3).
As corollary, for any constant $A$,
the set of all stable curves $X$ over $O_K$ with
$$
\frac{(\achern{1}{\omega_{X/O_K}, \Phi_{can}}^2)}{[K : \QQ]} \leq A
$$
is finite under the following equivalence (cf. Theorem~3.1).
For stable curves $X$ and $Y$ over $O_K$, $X$ is equivalent to $Y$
if $X \otimes_{O_K} O_{K'} \simeq Y \otimes_{O_K} O_{K'}$
for some finite extension field $K'$ of $K$.

In \S1, we will consider semistability of the kernel
of $H^0(C, L) \otimes \OO_C \to L$, which gives a generalization of
\cite{PR}. In \S2, an inequality of self-intersection and height will be
treated. Finally, \S3 is devoted to
finiteness of stable arithmetic surfaces with
bounded self-intersections of dualizing sheaves

We would like to thank Professor J.-B. Bost, S. Lang and
L. Szpiro for their helpful comments.

\head
1. Semistability of the kernel of $H^0(C, L) \otimes \OO_C \to L$
\endhead

Throughout this section, we will fix an algebraically closed
field $k$.
Let $C$ be a smooth projective curve over $k$.
For a non-zero torsion free sheaf $E$ on $C$, an average
degree $\mu(E)$ of $E$ is defined by $\mu(E) = \deg(E)/\rank E$.
$E$ is said to be stable (resp. semistable) if
$\mu(F) < \mu(E)$ (resp. $\mu(F) \leq \mu(E)$)
for all non-zero proper subsheaves $F$ of $E$.

Let $L$ a line bundle on $C$.
We set
$$
E(L) = \ker(H^0(C, L) \otimes \OO_C \to L)\quad\hbox{and}\quad
\Mov(L) = \image(H^0(C, L) \otimes \OO_C \to L).
$$
Clearly, $h^0(C, L) = h^0(C, \Mov(L))$.
If $h^0(C, L) \geq 2$, then $E(L) \not= 0$ and
$$
\mu(E(L)) = \frac{- \deg(\Mov(L))}{h^0(C, L) - 1}
= \frac{- \deg(\Mov(L))}{h^0(C, M(L)) - 1}.
$$
Moreover, $\mu(E(L)) \geq -2$ if and only if
$h^0(C, L) \geq \frac{1}{2} \deg(\Mov(L)) + 1$.
The main purpose of this section is to give a generalization
of A. Paranjape and S. Ramanan's result \cite{PR}.

\proclaim{Theorem 1.1}
Let $C$ be a smooth projective curve of genus $g \geq 1$ over $k$ and
$L$ a line bundle on $C$ such that $h^0(C, L) \geq 2$
and $\mu(E(L)) \geq -2$.
Then, $E(L)$ is semistable. Moreover, we have the following.
\roster
\item "(1)" If $\deg L \geq 2g + 1$, $E(L)$ is stable.

\item "(2)" In the case where $\deg L = 2g$ and $C$ is not hyperelliptic,
$E(L)$ is stable if and only if $h^0(C, L \otimes \omega_C^{-1}) = 0$.

\item "(3)" $E(\omega_C)$ is stable if and only if
$C$ is not hyperelliptic, where $\omega_C$ is the dualizing
sheaf of $C$ over $k$.
\endroster
\endproclaim

\demo{Proof}
First of all, we will prepare several lemmas.

\proclaim{Lemma 1.2}
Let $X$ be a $d$-dimensional projective variety over $k$ and $E$ a vector
bundle of rank $r$ on $X$. If $E$ is generated by global sections,
then there is a subvector space $V$ of $H^0(X, E)$ over $k$ such that
$\dim_k V = d + r$ and $V \otimes \OO_X \to E$ is surjective.
\endproclaim

\demo{Proof}
First, we consider a case where $r = 1$.
Let $\varphi : X \to \PP(H^0(X, E))$ the morphism induced
by the complete linear system $|E|$.
Since $\dim \varphi(X) \leq d$, there is a linear subspace
$T$ of $\PP(H^0(X, E))$ such that $\codim T = d + 1$ and
$\varphi(X) \cap T = \emptyset$.
Let $\{ s_0, s_1, \ldots, s_d \}$ be
a basis of $H^0(\PP^d, \OO_{\PP^d}(1))$
and $\pi : \varphi(X) \to \PP^{d}$ the morphism
induced by the projection $\PP(H^0(X, E)) \setminus T \to \PP^d$.
Here we consider a subvector space $V$ of $H^0(X, E)$
generated by
$$
(\pi \cdot \varphi)^*(s_0),(\pi \cdot \varphi)^*(s_1),
\ldots, (\pi \cdot \varphi)^*(s_d).
$$
Then, it is easy to see that $V$ is a desired vector
subspace.

Next, we consider a general case. Let $f : Y = \PP(E) \to X$
be the projective bundle of $E$ and $\OO_Y(1)$ the
tautological line bundle of $Y$.
Since $E$ is generated by global sections, so is $\OO_Y(1)$.
Thus, by the previous observation,
there is a subvector space $V$ of $H^0(Y, \OO_Y(1))$
such that $\dim V = d + r$ and $V \otimes \OO_Y \to \OO_Y(1)$
is surjective. Since $H^0(Y, \OO_Y(1)) \simeq H^0(X, E)$,
we can view $V$ as a subvector space of $H^0(X, E)$.
Pick up a point $x \in X$.
Let us consider the natural homomorphism $V \to H^0(Y_x, \OO_{Y_x}(1))$.
This is surjective because
$V \otimes \OO_{Y_x} \to \OO_{Y_x}(1)$ is surjective and
$W \otimes \OO_{Y_x} \to \OO_{Y_x}(1)$ is not surjective
for every proper subvector space $W$ of $H^0(Y_x, \OO_{Y_x}(1))$.
Therefore, $V \to E_x$ is surjective. Thus, we get our assertion.
\qed
\enddemo

\proclaim{Lemma 1.3}
Let $C$ be a smooth projective curve of genus $g$ over $k$ and
$W$ a non-zero vector bundle on $C$ such that $h^0(C, W) = 0$
and $W^*$ is generated by global sections.
Then, we have the following.
\roster
\item "(1)" If $h^1(C, (\det W)^{-1}) = 0$, then
$\deg W \leq - \rank W - g$.

\item "(2)" If $h^1(C, (\det W)^{-1}) \not= 0$, then
$\deg W \leq -2 \rank W$. Moreover, equality holds if and only if
either $W \simeq E(\omega_C)$, or $C$ is hyperelliptic and
$W \simeq E(\OO_C(m g^1_2))$ for some $1 \leq m \leq g-1$, where
$g^1_2$ is the hyperelliptic divisor.
\endroster
\endproclaim

\demo{Proof}
By Lemma~1.2, there is a surjective homomorphism
$\OO_C^{\rank W + 1} \to W^*$. Clearly, the kernel of it
is $\det W$. Thus, we get
$$
0 \to W \to \OO_C^{\rank W + 1} \to (\det W)^{-1} \to 0.
\tag 1.3.1
$$
Therefore, since $h^0(C, W) = 0$, we have
$$
    \rank W + 1 \leq h^0(C, (\det W)^{-1}).
\tag 1.3.2
$$

If $h^1(C, (\det W)^{-1}) = 0$, then, by Riemann-Roch theorem,
$$
h^0(C, (\det W)^{-1}) = - \deg W + 1 - g.
$$
Thus, by (1.3.2), we obtain (1).

If $h^1(C, (\det W)^{-1}) \not= 0$,
then $(\det W)^{-1}$ is special. Thus,
by Clifford's theorem (cf. Chapter IV, Theorem~5.4 of \cite{Ha}),
$$
h^0(C, (\det W)^{-1}) \leq \frac{-\deg W}{2} + 1.
\tag 1.3.3
$$
Therefore, by (1.3.2) and (1.3.3), we have $\deg W \leq -2 \rank W$.

If $\deg W = -2 \rank W$ holds, then we have
$$
    \rank W + 1 = h^0(C, (\det W)^{-1}) \quad\hbox{and}\quad
   h^0(C, (\det W)^{-1}) = \frac{-\deg W}{2} + 1.
$$
By equality conditions of Clifford's theorem, we have either
$\det W = \OO_C$, $\det W = \omega_C^{-1}$, or
$C$ is hyperelliptic and $\det W = \OO_C(-m g^1_2)$
for some $1 \leq m \leq g-1$.
The case $\det W = \OO_C$ is impossible because $\rank W = 0$ in this case.
If $\det W = \omega_C^{-1}$, then $\rank W = g -1$. Therefore,
by the exact sequence (1.3.1), we have $W \simeq E(\omega_C)$.
By the same way,
if $C$ is hyperelliptic and $\det W = \OO_C(-m g^1_2)$
for some $1 \leq m \leq g-1$, then $W \simeq E(\OO_C(m g^1_2))$.

Conversely, if $W \simeq E(\omega_C)$, or $W \simeq E(\OO_C(m g^1_2))$,
then it is easy to see that $\deg W = -2 \rank W$
\qed
\enddemo

\proclaim{Lemma 1.4}
Let $C$ be a smooth projective curve of genus $g$ over $k$ and
$L$ a line bundle on $C$. If $L$ is generated by
global sections and $h^0(C, L) \geq 2$, then
$$
h^0(C, E(L)^*) \geq h^0(C, L).
$$
Moreover, equality holds if and only if
$$
H^0(C, L) \otimes H^0(C, \omega_C) \to H^0(C, L \otimes \omega_C)
$$
is surjective. In particular, if $g \geq 2$ and $C$ is not hyperelliptic,
then $h^0(C, E(\omega_C)^*) = g$.
\endproclaim

\demo{Proof}
Let us consider the exact sequence:
$$
0 \to L^{-1} \to H^0(C, L)^* \otimes \OO_C \to E(L)^* \to 0.
$$
Thus, we have
$$
h^0(C, E(L)^*) \geq h^0(C, L).
$$
Further, equality holds if and only if
$$
H^1(C, L^{-1}) \to H^1(C, H^0(C, L)^* \otimes \OO_C)
$$
is injective. Considering Serre's duality, the injectivity is
equivalent to the surjectivity of
$$
H^0(C, H^0(C, L) \otimes \omega_C) \to H^0(C, L \otimes \omega_C).
$$
Thus, we have our lemma because
$H^0(C, H^0(C, L) \otimes \omega_C) \simeq
H^0(C, L) \otimes H^0(C, \omega_C).$
\qed
\enddemo

Let us start the proof of Theorem~1.1.
Let $W$ be a proper non-zero subvector bundle of $E(L)$. Then,
clearly, $h^0(C, W) = 0$ because $h^0(C, E(L)) = 0$.
Moreover, since $\Mov(L)$ is locally free, $E(L)^*$ is a quotient of
$H^0(C, L)^* \otimes \OO_C$. Thus so is $W^*$.
Therefore, $W^*$ is generated by global sections.
So we can apply Lemma~1.3.
If $h^1(C, (\det W)^{-1}) = 0$, then we have $\mu(W) < \mu(E(L))$
as follows.
$$
\mu(W) \leq -1 - \frac{g}{\rank W} < -1 - \frac{g}{\rank E(L)}
        = \frac{-\deg(\Mov(L)) - h^1(C, \Mov(L))}{\rank E(L)}
        \leq \mu(E(L)).
$$
If $h^1(C, (\det W)^{-1}) \not= 0$, then
$\mu(W) \leq -2$ by Lemma~1.3. Thus, $E(L)$ is semistable.

Next, we will consider stability of $E(L)$ for each case (1) -- (3).

\medskip
(1) In this case, $\mu(E(L)) > -2$. Thus, stability is trivial.

\medskip
(2) First, we assume that $E(L)$ is stable.
If $h^0(C, L \otimes \omega_C^{-1}) \not= 0$, then $\omega_C$
is a subsheaf of $L$. Thus, $E(\omega_C)$ is a subsheaf of $E(L)$.
On the other hand, $\mu(E(L)) = \mu(E(\omega_C)) = -2$. Thus, $E(L)$ is not
stable.
This is a contradiction.

Next we assume that $h^0(C, L \otimes \omega_C^{-1}) = 0$.
If $E(L)$ is not stable, by Lemma~1.3, there is a subbundle $W$
of $E(L)$ such that $W$ is isomorphic to $E(\omega_C)$.
Thus, we have an exact sequence:
$$
0 \to L \otimes \omega_C^{-1} \to E(L)^* \to E(\omega_C)^* \to 0.
$$
By Lemma~1.4, $h^0(C, E(L)^*) \geq g+1$ and $h^0(C, E(\omega_C)^*) = g$.
Thus, $h^0(C, L \otimes \omega_C^{-1}) \not= 0$.
This is a contradiction.

\medskip
(3) By Lemma~1.3, it is easy to see that if $C$ is not hyperelliptic,
then $E(\omega_C)$ is stable.
Here, we assume that $C$ is hyperelliptic.
Then, $\omega_C \simeq \OO_C((g-1)g^1_2)$.
Thus, $\omega_C$ has a subsheaf $\OO_C(g^1_2)$, which implies that
$E(\omega_C)$ has a subsheaf $E(g^1_2)$. On the other hand,
$\mu(E(\omega_C)) = \mu(E(g^1_2)) = -2$. Thus, $E(\omega_C)$
is not stable.
\qed
\enddemo

\head
2. Inequality of self-intersection and height
\endhead

Let $K$ be a number field and $O_K$ the ring of integers of $K$.
Let us consider a pair $(V, h)$ of
an $O_K$-module $V$ of finite rank and
Hermitian metric $h_{\sigma}$ on $V_{\sigma}$ for each
$\sigma \in K_{\infty}$.
We define $L^2$-degree
$\deg_{L^2}(V, h)$ of $(V, h)$ by
$$
  \deg_{L^2}(V, h) = \log\#\left( \frac{V}{O_K x_1 + \cdots + O_K x_t} \right)
   - \frac{1}{2} \sum_{\sigma \in K_{\infty}}
     \log \det(h_{\sigma}(x_i, x_j)),
$$
where $x_1, \ldots, x_t \in V$ and $\{ x_1, \ldots, x_t \}$ is a basis
of $V \otimes K$. Using the Hasse product formula,
it is easily checked that $\deg_{L^2}(V, h)$ does not depend on the choice of
$\{ x_1, \ldots, x_t \}$.

The purpose of this section is to prove the following theorem,
which is a variant of Theorem II in \cite{Bo} and
gives a refine result in special cases.

\proclaim{Theorem 2.1}
Let $f : X \to S = \Spec(O_K)$ be a regular
arithmetic surface of genus $g \geq 1$
and $L$ a line bundle on $X$ such that
$L$ is $f$-nef and $\deg(L_{\overline{\QQ}}) > 0$.
Let $V$ be a $O_K$-submodule of $H^0(X, L)$ and
$h$ a Hermitian metric of $V$ such that
the natural homomorphism
$V_{\overline{\QQ}} \otimes \OO_{X_{\overline{\QQ}}}
\to L_{\overline{\QQ}}$ is surjective.
Let $h_L$ the quotient metric of $L$ induced by $h$ via
the surjective homomorphism
$V_{\sigma} \otimes \OO_{X_{\sigma}} \to L_{\sigma}$
on each infinite fiber $X_{\sigma}$.
If $\ker(V_{\overline{\QQ}} \otimes \OO_{X_{\overline{\QQ}}}
\to L_{\overline{\QQ}})$ is semistable, then
$$
\frac{1}{2} \frac{(\achern{1}{L, h_L}^2)}{\deg(L_{\overline{\QQ}})} >
\frac{\deg_{L^2}(V, h)}{\rank V}.
$$
\endproclaim

\demo{Proof}
Let $Q$ be the image of
$V \otimes \OO_X \to L$ and
$S$ the kernel of
$V \otimes \OO_X \to L$.
Then, $Q$ is a torsion free sheaf of rank $1$ and
$S$ is a torsion free sheaf of rank $\rank V -1$.
Using the natural metric $f^*(h)$ of
$V \otimes \OO_X$,
we can give the quotient metric $h_Q$ to $Q$ and
the submetric $h_S$ to $S$.
Clearly, on each infinite fiber, $h_Q$ coincides with
$h_L$ of $L$.

Here, we calculate $\achern{1}{S, h_S}$ and
$\achern{2}{S, h_S}$ in terms of
$\achern{1}{V, h}$,
$\achern{1}{Q, h_Q}$ and the extension class of
$0 \to S \to V \otimes \OO_X \to Q \to 0$.
First of all, we get
$$
   \achern{1}{S, h_S} =
   f^*(\achern{1}{V, h}) - \achern{1}{Q, h_Q}.
\tag 2.1.1
$$
We set
$$
\rho =
\achern{2}{f^*(V, h)} -
\achern{2}{(S, h_S) \oplus (Q, h_Q)}.
\tag 2.1.2
$$
Then, by Proposition~7.3 of \cite{Mo},
we have that $\rho \geq 0$, and
$\rho = 0$ if and only if
the exact sequence:
$$
0 \to (S, h_{S}) \to
f^*(V, h)
\to (Q, h_Q) \to 0
\tag 2.1.3
$$
splits orthogonally on each infinite fiber.
It follows from (2.1.2) that
$$
\achern{2}{S, h_S} =
-(\rho + \achern{2}{Q, h_Q}) -
\deg(L_{\overline{\QQ}}) \deg_{L^2}(V, h) +
\achern{1}{Q, h_Q}^2.
\tag 2.1.4
$$
In our situation, (2.1.3) doesn't split on each infinite fiber
because a trivial bundle doesn't have an ample sub-line bundle
as its direct summand. So we get $\rho > 0$.
Moreover, since $\rank Q = 1$, we obtain
$\achern{2}{Q, h_Q} \geq 0$.
Hence, by (2.1.4),
$$
\achern{2}{S, h_S} < - \deg(L_{\overline{\QQ}}) \deg_{L^2}(V, h)
+ \achern{1}{Q, h_Q}^2,
\tag 2.1.5
$$

Since $S_{\overline{\QQ}}$ is semistable vector bundle,
by virtue of Corollary~8.9 in \cite{Mo}, we obtain
$$
(\rank V -2)\achern{1}{S, h_S}^2 \leq 2(\rank V -1) \achern{2}{S, h_S}.
\tag 2.1.6
$$
Combining (2.1.1), (2.1.5) and (2.1.6), we have
$$
\frac{1}{2} \frac{\achern{1}{Q, h_Q}^2}{\deg(L_{\overline{\QQ}})} >
\frac{\deg_{L^2}(V, h)}{\rank V }.
\tag 2.1.7
$$
On the other hand, since $Q \subseteq L$ and $Q \otimes K = L \otimes K$,
there is a vertical effective $1$-cycle $Z$ on $X$ such that
$$
\achern{1}{Q, h_Q} = \achern{1}{L, h_L} - Z.
$$
Therefore,
$$
\achern{1}{Q, h_Q}^2 = \achern{1}{L, h_L}^2 -
2 (L \cdot Z) + Z^2.
$$
Hence, since $L$ is $f$-nef and $Z^2 \leq 0$,
the above implies that
$$
\achern{1}{Q, h_Q}^2 \leq \achern{1}{L, h_L}^2.
\tag 2.1.8
$$
Thus, by (2.1.7) and (2.1.8),
we finally get our inequality.
\qed
\enddemo

By Theorem~1.1 and Theorem~2.1, we have

\proclaim{Corollary 2.2}
Let $f : X \to S = \Spec(O_K)$ be a regular
arithmetic surface of genus $g \geq 1$
and $L$ a line bundle on $X$ such that
$L$ is $f$-nef, $\deg(L_{\overline{\QQ}}) > 0$,
$L_{\overline{\QQ}}$ is generated by global sections
and $\rank H^0(X, L) \geq \frac{1}{2} \deg(L_{\overline{\QQ}}) + 1$.
Let $h$ be a Hermitian metric of $H^0(X, L)$ and
$h_L$ the quotient metric of $L$ induced by $h$ via
the surjective homomorphism
$H^0(X_{\sigma}, L_{\sigma}) \otimes \OO_{X_{\sigma}} \to L_{\sigma}$
on each infinite fiber.
Then we have
$$
\frac{1}{2} \frac{(\achern{1}{L, h_L}^2)}{\deg(L_{\overline{\QQ}})} >
\frac{\deg_{L^2}(H^0(X, L), h)}{\rank H^0(X, L)}.
$$
\endproclaim

Let $C$ be a compact Riemann surface of genus $g \geq 1$ and
$\Omega_C$ the sheaf of holomorphic $1$-forms on $C$.
The natural Hermitian metric $\langle \ , \ \rangle_{can}$ of
$H^0(C, \Omega_C)$ is defined by
$$
\langle \alpha , \beta \rangle_{can} =
\frac{\sqrt{-1}}{2} \int_C \alpha \wedge \bar{\beta}.
$$
Since $H^0(C, \Omega_C) \otimes \OO_C \to \Omega_C$ is surjective,
the Hermitian metric $\langle \ , \ \rangle_{can}$
induces the quotient Hermitian metric of
$\Omega_C$. We denote this metric by $\Phi_{can}$ and call it
the canonical metric of $\Omega_C$.
Let $\{ \omega_1, \ldots, \omega_g \}$ be an orthonormal basis
of $H^0(C, \Omega_C)$ with respect to $\langle \ , \ \rangle_{can}$.
Then, the K\"{a}hler metric $k_{can} = \Phi_{can}^{-1}$ is given by
$$
    \omega_1 \otimes \bar{\omega}_1 + \cdots +
    \omega_g \otimes \bar{\omega}_g.
$$

Let $K$ be a number field, $O_K$ the ring of integers of $K$
and $S = \Spec(O_K)$.
Let $f : X \to S$ be an arithmetic surface of the genus $g \geq 1$
with the invertible dualizing sheaf $\omega_{X/S}$.
We can give the above canonical metric to
$\omega_{X/S}$ on each infinite fiber.
By abuse of notation, we denote this metric by $\Phi_{can}$.

Let $A$ be an abelian variety of dimension $g$ over $K$ such that
$A$ has semi-stable reduction.
Let $\pi : N(A) \to S$ be the Neron model of $A$
and $\varepsilon : S \to N(A)$ the identity
of the group scheme $N(A)$.
Set
$\omega_{A/S} = \varepsilon^*(\det(\Omega_{N(A)/S}))$.
For each infinite place $\sigma \in K_{\infty}$,
we give a Hermitian metric $\langle \ , \ \rangle_{\sigma}$
of $\omega_{A/S}$ defined by
$$
\langle \alpha, \beta \rangle_{\sigma}
= \left( \frac{\sqrt{-1}}{2} \right)^g \int_{A_{\sigma}}
\alpha \wedge \bar{\beta}.
$$
The Faltings modular height $\Height_{Fal}(A)$ of $A$ is given by
$$
\Height_{Fal}(A) =
\frac{\deg_{L^2}(\omega_{A/S}, \langle , \rangle)}{[K : \QQ]}.
$$

The following corollary is the main result of this note
which is a refinement of the inequality (4.9) in \cite{Bo}.

\proclaim{Corollary 2.3}
Let $K$ be a number field, $O_K$ the ring of integers of $K$
and $S = \Spec(O_K)$.
Let $f : X \to S$ be a stable
arithmetic surface of genus $g \geq 2$.
Then we have
$$
\frac{(\achern{1}{\omega_{X/S}, \Phi_{can}}^2)}{[K : \QQ]} >
\frac{4(g-1)}{g} \Height_{Fal}(J(X_K)),
$$
where $J(X_K)$ is the Jacobian of $X_K$.
\endproclaim

\demo{Proof}
First of all, it is well known that
$$
\frac{\deg_{L^2}(H^0(X, \omega_{X/S}), \langle , \rangle_{nat})}
     {[K : \QQ]}
= \Height_{Fal}(J(X_K)).
$$
Let $\mu : Y \to X$ be a minimal resolution of $X$.
Then, $\omega_{Y/S}$ is $f$-nef,
$H^0(Y, \omega_{Y/S}) = H^0(X, \omega_{X/S})$, and
$(\achern{1}{\omega_{X/S}, \Phi_{can}}^2) =
(\achern{1}{\omega_{Y/S}, \Phi_{can}}^2)$.
Thus, our assertion follows from Corollary~2.2
\qed
\enddemo

\head
3. Finiteness of stable arithmetic surfaces with
bounded self-intersections of dualizing sheaves
\endhead

In this section, we will consider an application of the inequality of
Corollary~2.3.
Let $g$ an integer with $g \geq 2$. For a scheme $T$, we set
$$
\bar{M}^s_g(T) = \{ X \to T \mid \hbox{$X \to T$ is a stable
curve of genus $g$ with smooth generic fibers} \}.
$$
For $X, Y \in \bar{M}^s_g(\Spec(O_K))$, we define
an equivalence $X \sim Y$ by the following:
$$
X \sim Y \quad\underset\hbox{def}\to\Longleftrightarrow\quad
X \otimes_{O_K} O_{K'} \simeq Y \otimes_{O_K} O_{K'}
\quad\hbox{for some finite extension field $K'$ of $K$}.
$$
For a constant $A$, we denote by $B_g^{can}(K, A)$
a subset of $\bar{M}^s_g(\Spec(O_K))/\!\!\sim$
consisting of the classes of stable curves with
$$
\frac{(\achern{1}{\omega_{X/S}, \Phi_{can}}^2)}{[K : \QQ]} \leq A.
$$
Then, we have

\proclaim{Theorem~3.1}
If $g \geq 2$, then
$B_g^{can}(K, A)$ is finite for any number field $K$ and
any constant $A$.
\endproclaim

\demo{Proof}
By Corollary~2.3 and the finiteness property of
the Faltings modular height (cf. \cite{Fa}), it is sufficient to prove
the following lemma.
\enddemo

\proclaim{Lemma~3.2}
Let $K$ be a number field and $O_K$ the ring of integers of $K$.
Let $X$ and $X'$ be stable curves over $O_K$ of genus $g \geq 2$.
If $X_K$ is isomorphic to $X'_K$ over $K$, then
this isomorphism extends to an isomorphism over $O_K$.
\endproclaim

\demo{Proof}
Since $X_K$ is isomorphic to $X'_K$ over $K$, there is
a rational map $\phi : X \dashrightarrow X'$ over $O_K$.
Let $\mu : Y \to X$ be a minimal succession of blowing-ups
such that $\phi \cdot \mu$ induces a morphism $\mu' : Y \to X'$.
There are effective divisors $Z$ and $Z'$ on $Y$
such that
$$
\omega_{Y/O_K} = \mu^*(\omega_{X/O_K}) \otimes \OO_Y(Z)
= {\mu'}^*(\omega_{X'/O_K}) \otimes \OO_Y(Z').
$$
Clearly, $\Supp(Z)$ is contracted by $\mu$
and $\Supp(Z')$ is contracted by $\mu'$.
Moreover, by the minimality of $\mu$, $Z$ and $Z'$ has
no common components.
Let us consider
$$
  (\mu^*(\omega_{X/O_K}) \otimes \OO_Y(Z) \cdot Z) =
  ({\mu'}^*(\omega_{X'/O_K}) \otimes \OO_Y(Z') \cdot Z).
$$
If $Z \not= 0$, then the left hand side of the above is negative.
But the right hand side is non-negative. Therefore, $Z = 0$.
By the same way, $Z' = 0$. Thus, we have
$$
\omega_{Y/O_K} = \mu^*(\omega_{X/O_K}) = {\mu'}^*(\omega_{X'/O_K}).
$$
Here, we assume that $\mu$ is not an isomorphism.
Then, there is a curve $C$ on $Y$ such that $C$ is contracted by $\mu$,
but is not contracted by $\mu'$.
Then,
$(\mu^*(\omega_{X/O_K}) \cdot C) = 0$, but
$({\mu'}^*(\omega_{X'/O_K}) \cdot C) > 0$.
This is a contradiction.
Thus, the rational map $\phi : X \dashrightarrow X'$ is actually
a morphism and $\phi^*(\omega_{X'/O_K}) = \omega_{X/O_K}$.
Hence, $\phi$ is finite
because $\omega_{X/O_K}$ and $\omega_{X'/O_K}$ are ample.
Moreover, $X$ and $X'$ are normal. Therefore, by Zariski main theorem,
$\phi$ is an isomorphism.
\qed
\enddemo

Next, we give a variant of Theorem~3.1.
Let $C$ be a compact Riemann surface of genus $g \geq 1$ and
$\{ \omega_1, \cdots, \omega_g \}$ an orthonormal basis of $H^0(C, \Omega_C)$
with respect to $\langle,\rangle_{can}$.
We set the normalized K\"{a}hler form $\mu_{can}$ as follows.
$$
   \mu_{can} = \frac{\sqrt{-1}}{2 g} \sum_{i=1}^g
   \omega_i \wedge \bar{\omega}_i.
$$

We can give another metric $\Phi_{Ar}$
of $\Omega_C$ by the following way, which
is called the Arakelov metric. Let $\Delta$
be the diagonal of $C \times C$,
$p : C \times C \to C$ the first projection, and
$q : C \times C \to C$ the second projection.
Let $h_{\Delta}$ be the Einstein-Hermitian metric of
$\OO_{C \times C}(\Delta)$ with respect to
$p^* \mu_{can} + q^* \mu_{can}$ such that
$$
\int_{C \times C} \log(h_{\Delta}(1, 1))
(p^* \mu_{can} + q^* \mu_{can})^2 = 0,
$$
where $1$ is the canonical section of $\OO_{C \times C}(\Delta)$.
Since $\Omega_C$ is canonically isomorphic to $\OO_{\Delta}(-\Delta)$,
$\rest{(h_{\Delta}^{-1})}{\Delta}$
induces the metric $\Phi_{Ar}$ on $\Omega_C$.
It is well known that $c_1(\Omega_C, \Phi_{Ar}) = 2(g-1)\mu_{can}$.
Here we set
$$
\nu_{Ar}(C) = \int_C \log\left(\frac{\Phi_{Ar}}{\Phi_{can}}\right) \mu_{can}.
$$
Then, we have the following lemma.

\proclaim{Lemma 3.3}
Let $K$ be a number field, $O_K$ the ring of integers of $K$
and $S = \Spec(O_K)$.
Let $f : X \to S$ be a stable arithmetic surface of genus $g \geq 1$.
Then,
$$
(\achern{1}{\omega_{X/S}, \Phi_{can}}^2) =
(\achern{1}{\omega_{X/S}, \Phi_{Ar}}^2) + 2(g-1)\nu_{Ar}(X/S),
$$
where $\nu_{Ar}(X/S) = \sum_{\sigma \in K(\CC)} \nu_{Ar}(X_{\sigma})$.
\endproclaim

\demo{Proof}
We set $\rho = \Phi_{Ar}/\Phi_{can}$. Then
$$
(\achern{1}{\omega_{X/S}, \Phi_{can}}^2) =
(\achern{1}{\omega_{X/S}, \Phi_{Ar}}^2) + \sum_{\sigma \in K(\CC)}
\int_{C_{\sigma}} \log(\rho_{\sigma})
c_1(\omega_{X_{\sigma}}, (\Phi_{Ar})_{\sigma}).
$$
Thus we get our lemma because
$c_1(\omega_{X_{\sigma}}, (\Phi_{Ar})_{\sigma}) =
2(g-1)(\mu_{can})_{\sigma}$.
\qed
\enddemo

For constants $A_1$ and $A_2$,
we denote by $B_g^{Ar}(K, A_1, A_2)$
a subset of $\bar{M}^s_g(\Spec(O_K))/\!\!\sim$
consisting of the classes of stable curves with
$$
\frac{(\achern{1}{\omega_{X/S}, \Phi_{Ar}}^2)}{[K : \QQ]} \leq A_1
\quad\hbox{and}\quad
\frac{\nu_{Ar}(X/S)}{[K : \QQ]} \leq A_2.
$$
Thus, Theorem~3.1 and Lemma~3.3 imply

\proclaim{Corollary~3.4}
For any number field $K$ and
any constants $A_1$ and $A_2$, $B_g^{Ar}(K, A_1, A_2)$
is finite.
\endproclaim

\definition{Remark~3.5}
Optimistically, we can guess that
$$
B_g^{Ar}(K, A) = \left\{ X \to S \in \bar{M}^s_g(\Spec(O_K)) \ \left|
\ \frac{(\achern{1}{\omega_{X/S}, \Phi_{Ar}}^2)}{[K : \QQ]} \leq A
\right\}\right/\!\!\sim
$$ is finite for any constant $A$.
Indeed, as L. Szpiro pointed out in his letter,
$B_g^{Ar}(K, 0)$ is finite.
For, let $f : X \to S$ be a stable arithmetic surface of genus $g \geq 2$.
Then, by \cite{Zh}, if $f$ is not smooth, then
$(\achern{1}{\omega_{X/S}, \Phi_{Ar}}^2) > 0$.
Thus, by Shafarevich Conjecture which was proved by Faltings \cite{Fa},
$B_g^{Ar}(K, 0)$ is finite.
\enddefinition


\enddocument